\documentclass[12pt]{article}
\usepackage{latexsym}
\usepackage{amsfonts}
\usepackage{amsmath} 
\usepackage{amssymb}
\usepackage{amsthm}
\usepackage{parskip}
\usepackage{enumerate}
 \usepackage{hyperref} 

\newcommand{\cla}{{\mathcal A}}
\newcommand{\clb}{{\mathcal B}}
\newcommand{\clc}{{\mathcal C}}
\newcommand{\cld}{{\mathcal D}}
\newcommand{\cle}{{\mathcal E}}
\newcommand{\clf}{{\mathcal F}}
\newcommand{\clg}{{\mathcal G}}
\newcommand{\clh}{{\mathcal H}}

\newcommand{\clk}{{\mathcal K}}

\newcommand{\cln}{{\mathcal N}}

\newcommand{\clp}{{\mathcal P}}

\newcommand{\Eb}{{\mbox{\large ${\sf E}$}}}
\newcommand{\Pb}{{{\sf P}}}

\renewcommand{\L}{{\mathbb L}}

\newcommand{\W}{{\mathbb W}}

\newcommand{\Ind}{{\mathbf 1}}
\newcommand{\tca}{\phi}
\newcommand{\tcb}{\psi}
\newtheorem{lemma}{{\bf Lemma}}
\newtheorem{remark}[lemma]{{\bf Remark}}
\newtheorem{corollary}[lemma]{{\bf Corollary}}
\newtheorem{example}[lemma]{{\bf Example}}
   \newtheorem{definition}[lemma]{{\bf Definition}}
 \newtheorem{theorem}[lemma]{{\bf Theorem}}
 \newtheorem{proposition}[lemma]{{\bf Proposition}}
\begin{document} 

\centerline{\bf \large Stochastic Integrals and Two Filtrations}
\vskip 7mm

\centerline{\large Rajeeva L Karandikar}
\vskip 1mm
\centerline{and}
\vskip 1mm
\centerline{\large B. V. Rao}
\vskip 7mm

\centerline{\em Chennai Mathematical Institute}
\centerline{\em H1 Sipcot IT Park, Siruseri, TN 603103, India.} 

%\cortext[cor2]{Corresponding author}
\begin{abstract}
{In the definition of the stochastic integral, apart from the integrand and the integrator, there is an underlying filtration that plays a role. Thus, it is natural to ask: {\it Does the stochastic integral depend upon the filtration?} In other words, if we have two filtrations, $(\clf_\centerdot)$ and $(\clg_\centerdot)$, a process $X$ that is semimartingale under both the filtrations and a process $f$ that is predictable for both the filtrations, then are the two stochastic integrals - $Y=\int f\,dX$, with filtration $(\clf_\centerdot)$ and $Z=\int f\,dX$, with filtration $(\clg_\centerdot)$ the same?

 When $f$ is left continuous with right limits, then the answer is yes. When one filtration is an enlargement of the other, the two integrals are equal if $f$ is bounded but this may not be the case when $f$ is unbounded. 
 
 We discuss this and give sufficient conditions under which the two integrals are equal.}
\end{abstract}

\vfill

\hrule

\vspace*{3mm}

\noindent
{\bf Keywords}: Stochastic Integral,  Two filtrations,   Semimartingale\\
{\em email: rlk@cmi.ac.in, \hskip 7mm bvrao@cmi.ac.in}

\newpage

\section{Introduction}

Let $Y=\int f\,dX$ where $X$ is a semimartingale and $f$ is a predictable process. There is in the background, a filtration $ (\clf_t)_{t\ge 0}$, with $X$ being a semimartingale w.r.t.\! the filtration $(\clf_\centerdot)$ and $f$ is predictable for this filtration. A natural question is: {\em Does the integral $Y$ depend upon the filtration $(\clf_\centerdot)$?} In other words, if there is another filtration $ (\clg_t)_{t\ge 0}$ such that $X$ is a semimartingale w.r.t.\! the filtration $(\clg_\centerdot)$ and $f$ is predictable w.r.t.\! the filtration $(\clg_\centerdot)$ and writing as $Z$ the stochastic integral of $f$ w.r.t.\! $X$ with the underlying filtration being taken as $(\clg_\centerdot)$, the question is: are $Y$ and $Z$ equal?

The answer is known only in some specific cases- when one filtration is enlargement of the other ({\em i.e.} $\clf_t\subseteq \clg_t$ for all $t$), and $f$ is bounded then $Y$ and $Z$ are equal while this may not be so if $f$ is unbounded. In the general case, (when the two filtrations may not be comparable), it is known that when $f$ is bounded, $Y-Z$ is a continuous process and the quadratic variation of $Y-Z$ is zero. Further, if the jumps of $X$ are summable, then it is known that $Y-Z$ is a process with finite variation paths. 

We will give some conditions under which $Y$ and $Z$ are equal.

\section{Notations and Preliminaries}
Let $(\Omega, \clf, \Pb)$ be a complete probability space. 
Let $\cln=\{A\in\clf\,:\,\Pb(A)=0\}$ be the class of $\Pb$-null sets. For a collection of random variables $\{U_\alpha:\alpha\in\Delta\}$, $\sigma(U_\alpha:\alpha\in\Delta)$ will denote the smallest $\sigma$-filed containing $\cln$ with respect to which $\{U_\alpha:\alpha\in\Delta\}$ are measurable. Likewise, all filtrations $ (\clf_t)_{t\ge 0}$ we consider satisfy $\cln\subseteq\clf_0$. 
For definitions and classical results on stochastic integration, see Jacod \cite{J78}, Karandikar-Rao \cite{rlkbv2} or Protter \cite{Protter}. 

Since we will be working with two filtrations, we will write all notions such as predictable, martingale, local martingale, semimartingale with the underlying filtration as prefix: Thus, we will say that $f$ is $(\clf_\centerdot)$- predictable for a process $f$ that is predictable w.r.t.\! the filtration $(\clf_\centerdot)$, $X$ is a $(\clf_\centerdot)$- local martingale will mean that $X$ is a local martingale w.r.t.\! the filtration $(\clf_\centerdot)$ and so on. We will write
\[(\clf_\centerdot)\text{-}\!\!\int f\,dX\]
to denote the stochastic integral of a process $f$ w.r.t.\! $X$, when $f$ is $(\clf_\centerdot)$- predictable and $X$ is $(\clf_\centerdot)$- semimartingale and the stochastic integral is defined with the filtration $(\clf_\centerdot)$- in the picture - for example, defined for $(\clf_\centerdot)$- predictable simple processes and extended by continuity to suitable class of processes (that include bounded $(\clf_\centerdot)$- predictable processes). Let us denote by $\clp(\clf_\centerdot)$ - the predictable $\sigma$ field for the filtration $(\clf_\centerdot)$, $\W(\clf_\centerdot)$ the class of $\clp(\clf_\centerdot)$ measurable processes.  Let $\W_b(\clf_\centerdot)$ and $\W_l(\clf_\centerdot)$ denote the class of  bounded processes and locally bounded processes in $\W(\clf_\centerdot)$ respectively.

For an r.c.l.l.\! $(\clf_\centerdot)$- adapted process $H$, the process $H^-$ defined by
\[H^-(t,\omega)=\begin{cases}\lim_{u<t,\;\;u\rightarrow t}H(u,\omega) & \text{ if }t>0\\
0 &\text{ if }t=0\end{cases}\]
is $(\clf_\centerdot)$- adapted. Further, $H^-$ is left continuous and thus predictable. Further $H^-$ is locally bounded and $\int H^-\,dX$ is defined for every $(\clf_\centerdot)$- semimartingale $X$.

The quadratic variation of an $(\clf_\centerdot)$- semimartingale $X$ will be denoted by $[X,X]$. As noted in Karandikar-Rao \cite{rlkbv, rlkbv2}, the quadratic variation $[X,X]$ of a semimartingale $X$ can be defined pathwise and hence does not depend upon the underlying filtration. Likewise, the continuous part of the increasing process $[X,X]$, denoted by $[X,X]^c$ is defined by
\[[X,X]^c_t=[X,X]_t-\sum_{0\le s\le t}((\Delta X)_s)^2\]
and again is defined pathwise and thus $[X,X]^c$ does not depend upon the filtration.

Given an $(\clf_\centerdot)$- semimartingale $X$, there exists a unique continuous $(\clf_\centerdot)$- local martingale $M$ such that 
\[[X-M,U]=0\;\;\;\forall\text{ continuous $(\clf_\centerdot)$- local martingales }U,\]
see Karandikar-Rao \cite{rlkbv2}, Theorem 5.64.
$M$ is the continuous (local) martingale part of $X$ and is denoted by $X^c$, but since this depends upon the filtration, in this article we will denote this as
\[M=C(X,(\clf_\centerdot)).\]
It can be checked that
\begin{equation}\label{c1}
[X,X]^c_t=[M,M]_t\;\;\;\text{ where } M=C(X,(\clf_\centerdot)).\end{equation}
Also, if $Y=X+V$, where $V$ is a process whose paths have finite variation, then 
\[C(X,(\clf_\centerdot))=C(Y,(\clf_\centerdot)).\]
We will denote by $\L(X,(\clf_\centerdot))$ the class of $(\clf_\centerdot)$- predictable processes such that the integral $(\clf_\centerdot)\text{-}\!\!\int f\,dX$ is defined. Usually, $\L(X,(\clf_\centerdot))$ is described in terms of the decomposition of $X$ : $X=M+A$ where $M$ is an $(\clf_\centerdot)$- local martingale and $A$ is a process with finite variation paths. An equivalent way of describing this class, which plays an important role in Proposition \ref{pr2} below, is taken from Karandikar-Rao \cite{rlkbv2}.

The class of $X$-integrable processes - $\L(X,(\clf_\centerdot))$ consists of $(\clf_\centerdot)$- predictable processes $f$ such that
\begin{equation}\label{eq1}
h^n \text{ bounded }(\clf_\centerdot)\text{- predictable}, |h^n|\le |f|, \;h^n\rightarrow 0\text{ pointwise}\end{equation}
 implies that 
\[(\clf_\centerdot)\text{-}\!\!\int h^n\,dX\rightarrow 0 \text{ in {\em ucp} topology.}\] 
Further, for $f\in \L(X,(\clf_\centerdot))$,  \begin{equation}\label{h1}
\int f\Ind_{\{|f|\le n\}}\,dX\rightarrow \int f\,dX \text{ in {\em ucp} topology as $n\rightarrow\infty$}.\end{equation}

It can be seen that a locally bounded predictable process $f$ is  $X$-integrable  for every semimartinagle $X$, {\em i.e.} $f\in  \L(X,(\clf_\centerdot))$.

Here, convergence of processes $Z^n$ to $Z$ in {\em ucp} topology means
\[\lim_{n\rightarrow \infty}\Pb(\sup_{0\le t\le T}|Z^n_t-Z_t|>\epsilon)=0\;\;\;\;\forall T<\infty,\;\;\forall \epsilon>0.\] 
\section{A preliminary observation}

Let $X$ be a process with r.c.l.l.\! paths such that $X$ is a $(\clf_\centerdot)$- semimartingale as well as a $(\clg_\centerdot)$- semimartingale. The question we are considering is:\\ 
Let $f\in\W(\clf_\centerdot)\cap \W(\clg_\centerdot)$. Under what conditions on $f,\,X,\,(\clf_\centerdot),\,(\clg_\centerdot)$- does
\begin{equation}\label{qq0}
(\clf_\centerdot)\text{-}\!\!\int f\,dX=(\clg_\centerdot)\text{-}\!\!\int f\,dX\,?
\end{equation} 

More precisely, let $ f\in \W(\clf_\centerdot)\cap \W(\clg_\centerdot)$. 
\begin{enumerate}
\item If further, $f$ is bounded, is \eqref{qq0} true?
\item If $f\in \L(X, (\clf_\centerdot))$ then can we conclude that $f\in \L(X, (\clg_\centerdot))$ and then is \eqref{qq0} true?
\item If $f\in \L(X, (\clf_\centerdot))$ and $f\in \L(X, (\clg_\centerdot))$ then can we conclude that \eqref{qq0} is true?
\end{enumerate}
For a large class of integrands, the desired conclusion \eqref{qq0} is true as we observe first. This is a direct consequence of the {\em pathwise integration formula}: See Bichteler \cite{Bichteler}, Karandikar \cite{rlkthesis, rlk1, rlk5, rlk3, rlk4}. 
\vskip 5mm

\begin{proposition}\label{pr1}
Let $U$ be an r.c.l.l.\! process such that $U$ is $(\clf_\centerdot)$- adapted as well as $(\clg_\centerdot)$- adapted. Then 
\begin{equation}\label{qq1}
(\clf_\centerdot)\text{-}\!\!\int U^-\,dX=(\clg_\centerdot)\text{-}\!\!\int U^-\,dX\;\;\;
\end{equation} 
\end{proposition}
\begin{proof}
For each fixed $n$, define $\{\sigma^n_i :\;i\ge 0\}$ inductively with $\sigma^n_0=0$ and
\[\sigma^n_{i+1}=\inf\{ t> \sigma^n_i : |U_t-U_{\sigma^n_i }|\ge 2^{-n}\mbox{ or } |U_{t-}-U_{\sigma^n_i }|\ge 2^{-n} \}.
\]
 For all $n,i$, $\sigma^n_i$ is an $(\clf_\centerdot)$- stopping time as well as an $(\clg_\centerdot)$- stopping time. 
 Let
 \[Z^n_t= \sum_{j=0}^\infty U_{t\wedge\sigma^n_j} (X_{t\wedge\sigma^n_{j+1}}-X_{t\wedge\sigma^n_j} ).\] 
 Then (see Karandikar-Rao \cite{rlkbv2}, Theorem 6.2)
 \[(\clf_\centerdot)\text{-}\!\!\int U^-\,dX=\lim_{n\rightarrow \infty} Z^n \text{ in the {\em 
ucp} metric}\]
and also
 \[(\clg_\centerdot)\text{-}\!\!\int U^-\,dX=\lim_{n\rightarrow \infty} Z^n \text{ in the {\em 
ucp} metric}.\]
Thus \eqref{qq1} holds.
\end{proof}

\section{Case of Nested Filtrations} 
In this section, we consider the case when
\begin{equation}\label{z1x}
\clf_t\subseteq \clg_t\;\; \forall t.
\end{equation}

\vskip -4mm
Our first observation is:
\vskip 2mm

\begin{proposition}\label{pr2}
Suppose $ (\clf_\centerdot)$ and $ (\clg_\centerdot)$ satisfy \eqref{z1x}.
 \begin{enumerate}[(i)]
\item Let $f$ be locally bounded $(\clf_\centerdot)$ predictable process. Then
\begin{equation}\label{qq3}
(\clf_\centerdot)\text{-}\!\!\int f\,dX=(\clg_\centerdot)\text{-}\!\!\int f\,dX.
\end{equation} 
\item Let  $f\in\L(X, (\clg_\centerdot))$. Then $f\in\L(X, (\clf_\centerdot))$ and \eqref{qq3} is true.
\end{enumerate}
\end{proposition}
\begin{proof}
For $(i)$, note that \eqref{qq3} holds for simple $(\clf_\centerdot)$- predictable processes and hence by monotone class theorem, \eqref{qq3} holds for all bounded $(\clf_\centerdot)$- predictable processes. This is Theorem VIII. 13 in Dellacherie - Meyer \cite{delmey} . Thus \eqref{qq3} holds for bounded $f\in\W$.  For $f\in\W_l$, writing $f^n=f\Ind_{\{|f|\le n\}}$ and using \eqref{h1} holds for both the filtrations, we conclude that \eqref{qq3} holds for $f$.

 For $(ii)$ let $h^n$ be $(\clf_\centerdot)$- predictable bounded processes such that $|h^n|\le |f|$ and $h^n\rightarrow 0$ pointwise. Since $\clf_t\subseteq \clg_t$ for all $t$, it follows that $h^n$ are $(\clg_\centerdot)$- predictable. Using $f\in\L(X, (\clg_\centerdot))$, it follows that 
\[(\clg_\centerdot)\text{-}\!\!\int h^n\,dX\rightarrow 0 \;\;\text{ in {\em ucp} as }n\rightarrow \infty.\]
Moreover, for each $n$, $h^n$ is bounded and thus by part (i), the $(\clg_\centerdot)$ and $(\clf_\centerdot)$ integrals of $h^n$ (w.r.t.\! $X$) are identical. Thus
\[(\clf_\centerdot)\text{-}\!\!\int h^n\,dX\rightarrow 0 \;\text{ in {\em ucp} as }n\rightarrow \infty\]
and hence $f\in \L(X, (\clf_\centerdot))$. 
Writing $f^n=f\Ind_{\{|f|\le n\}}$, we have
\[(\clf_\centerdot)\text{-}\!\!\int f^n\,dX\rightarrow (\clf_\centerdot)\text{-}\!\!\int f\,dX \text{ in {\em ucp} as $n\rightarrow\infty$}\]
as well as 
\[(\clg_\centerdot)\text{-}\!\!\int f^n\,dX\rightarrow (\clg_\centerdot)\text{-}\!\!\int f\,dX \text{ in {\em ucp} as $n\rightarrow\infty$}.\]
Since $f^n$ are bounded, the $(\clg_\centerdot)$ and $(\clf_\centerdot)$ integrals of $f^n$ (w.r.t.\! $X$) are identical by part (i) and thus it follows that \eqref{qq3} holds for $f$. 
\end{proof}

\begin{remark} \label{rm2}
Suppose $\clf_t\subseteq \clg_t$ for all $t$. The condition $f\in L(X, (\clf_\centerdot))$ do not imply that $f\in L(X, (\clg_\centerdot))$. 
Let us recall example 5.77 from Karandikar-Rao \cite{rlkbv2}. 
Let $\{\xi^{k,m}:\,1\le k\le 2^{m-1}, m\ge 1\}$ be a family of independent identically distributed random variables with 
\[\Pb(\xi^{1,1}=1)=\Pb(\xi^{1,1}=-1)=0.5\]
and let $a^{k,m}=\frac{2k-1}{2^m}$. Let 
$\clf_t=\sigma\{\xi^{k,m}: \,a^{k,m}\le t\},$ 
\[A_t=\sum_{m=1}^\infty\sum_{k=1}^{2^{m-1}}\frac{1}{2^{2m}}\xi^{k,m}\Ind_{[a^{k,m},\infty)}(t)\]
and $f:[0,\infty)\mapsto [0,\infty)$ be defined by
\[f(a^{k,m})=2^m\]
with $f(t)=0$ otherwise.

It is shown in \cite{rlkbv2} that $\int f\,dA$ exists as a martingale integral, but does not exist as Riemann-Stieltjes integral. Now take $\clg_t=\clf_\infty$. Then it can be seen that the only $(\clg_\centerdot)$- local martingales are constants and thus $(\clg_\centerdot)\text{-}\!\!\int h\,dA$ for any $h$ is just the Riemann-Stieltjes integral. Thus $(\clg_\centerdot)\text{-}\!\!\int f\,dA$ does not exist while $(\clf_\centerdot)\text{-}\!\!\int f\,dA$ exists. 

%Also see \cite{Protter} p.361 for another example.
\end{remark}
\begin{example}\label{ex10}We now give the classical example, due to Ito \cite{Ito}: Let $W$ be a Brownian motion on some complete probability space and let $\clf_t=\sigma(W_s:0\le s\le t)$. Let $0<T<\infty$ be fixed and let  $(\clg_\centerdot)$ denote the filtration with $\clg_t=\sigma(W_s:0\le s\le t, \;W_{T})$ for $t\ge 0$. Ito showed that $W$ is not a $(\clg_\centerdot)$- martingale, but is a $(\clg_\centerdot)$- semimartingale with 
\[M_t=W_t-\int_0^{t\wedge T}\frac{W_T-W_s}{T-s}ds\]
being a martingale (indeed a Brownian motion) w.r.t. the filtration $(\clg_\centerdot)$. It has been shown (e.g. see p.361, \cite{Protter} ) that in this case $L(X, (\clf_\centerdot))$ is not a subset of $L(X, (\clg_\centerdot))$
\end{example}

We will now give a generic expansion of a filtration $(\clf_\centerdot)$ 
to $(\clg_\centerdot)$, where every $X$ that is a $(\clf_\centerdot)$ semimartingale is also a a $(\clg_\centerdot)$ semimartingale and 
\[L(X, (\clf_\centerdot))\subseteq L(X, (\clg_\centerdot)).\]
\vskip 3mm

\begin{theorem}\label{th4}
Let $\{A_k:\,k\ge 1\}$ be a partition (consisting of measurable sets) of $\Omega$. Let  $(\clf_\centerdot)$ be a filtration and let 
\begin{equation}\label{az1}
\clg_t=\sigma(\clf_t\cup \{A_k:\,k\ge 1\}) \;\;\;t\ge 0.
\end{equation}
Then if $X$ is a $(\clf_\centerdot)$- semimartingale, then $X$ is also a $(\clg_\centerdot)$- semimartingale. Further, in this case,
\[f\in L(X, (\clf_\centerdot)) \;\text{ implies }\; f\in L(X, (\clg_\centerdot))\]
and then \eqref{qq3} holds.
\end{theorem}
\begin{proof}
The first part is known as Jacod's countable expansion (see \cite{Protter} p.53). First let us observe that $g$ is $(\clg_\centerdot)$ - predictable if and only if it admits a representation
\begin{equation}\label{az0} 
g=\sum_{k=1}^\infty h^k\Ind_{A_k}
\end{equation}
where $h^k$ are $(\clf_\centerdot)$- predictable, and if $g$ is bounded by $M$ then $h^k$ can also be chosen to be bounded by $M$. This is easily verified for $(\clg_\centerdot)$ simple processes $g$. Let $g_n$ be processes such that
\begin{equation}\label{az19} 
g_n=\sum_{k=1}^\infty h^k_n\Ind_{A_k}
\end{equation}
where $g_n$ are uniformly bounded by $M$, $h^k_n$ for $k\ge 1$, $n\ge 1$ are $(\clf_\centerdot)$ - predictable also bounded by $M$, with $g_n$ converging pointwise to $g$. Let 
\[h^k=\limsup_{n\rightarrow\infty}h^k_n\]
Then $h^k, k\ge 1$ and $g$ would satisfy \eqref{az0} and thus the class of processes that admit a representation as in   \eqref{az0} is closed under bounded pointwise convergence and contains $(\clg_\centerdot)$ simple processes and thus equals the class of  bounded $(\clg_\centerdot)$ - predictable processes. Next we observe that  that if $g, \{h^k\}$ satisfy   
\eqref{az0}, then we have
\begin{equation}\label{az2} 
(\clg_\centerdot)\text{-}\!\!\int g\,dX=\sum_{k=1}^\infty \Ind_{A_k}\Bigl[(\clf_\centerdot)\text{-}\!\!\int h^k\,dX\Bigr].
\end{equation}

This can be verified for simple processes and then using the observation given above along with the dominated convergence theorem for stochastic integrals (see \cite{rlkbv2}, p.93), we conclude that \eqref{az2} holds for all  bounded $(\clg_\centerdot)$- predictable processes $g$ with $h^k$ as in \eqref{az0}.

Now let $f$ be a $(\clf_\centerdot)$- predictable processes such that $f\in\L(X,(\clf_\centerdot))$. To show that $f\in\L(X,(\clg_\centerdot))$, suffices to prove that if $g_n$ are bounded $(\clg_\centerdot)$ predictable processes, $|g_n|\le |f|$ for all $n$, and $g_n$ converges to 0 pointwise then $(\clg_\centerdot)\text{-}\!\!\int g_n\,dX$ converges to 0 (in ucp). For this, get $\{h^k_n\}$, for $n\ge 1, k\ge 1$ bounded $(\clf_\centerdot)$ predictable processes such that \eqref{az19} holds for all $n$. Then as noted above, we have
\begin{equation}\label{az22} 
(\clg_\centerdot)\text{-}\!\!\int g_n\,dX=\sum_{k=1}^\infty \Ind_{A_k}\Bigl[(\clf_\centerdot)\text{-}\!\!\int h^k_n\,dX\Bigr].
\end{equation}
 Let
\[B_k=\{\omega:\,\limsup_{n\rightarrow\infty}|h^k_n(\omega)|\neq 0\}.\]
Since $g_n$ converges to 0 pointwise, it follows that $A_k\cap B_k=\phi$. Replacing $h^k_n$ by $\max\{\min\{h^k_n,f\},f\}\Ind_{B^c_k}$, we have that \eqref{az19} holds,  $|h^k_n|\le |f|$ and for each $k$, $|h^k_n|$ converges to 0 pointwise as $n\rightarrow\infty$. Using $f\in\L(X,(\clf_\centerdot))$, it follows that
\[(\clf_\centerdot)\text{-}\!\!\int h^k_n\,dX\rightarrow 0\text{ in }ucp.\]
Now using \eqref{az2} and \eqref{az22}, we conclude that $(\clf_\centerdot)\text{-}\!\!\int g_n\,dX$ converges to 0 in ucp. This completes the proof. 
\end{proof}
\begin{remark}
The proof given above also includes the proof of Jacod's countable expansion theorem, namely that $X$ is a $(\clg_\centerdot)$- stochastic integrator and hence a $(\clg_\centerdot)$- semimartingale. See \cite{rlkbv2}.
\end{remark} 

The next remark follows from Proposition \ref{pr2} and Theorem \ref{th4}.\vskip 3mm

\begin{remark} \label{rem77}
Let  $(\clf_\centerdot)$, $(\clg_\centerdot)$ and $X$ be as in Theorem 
\ref{th4}. Let $(\clh_\centerdot)$ be a filtration such that
\[\clf_t\subseteq \clh_t\subseteq \clg_t\;\;\forall t\ge 0.\]
Then if $f\in \L(X,(\clf_\centerdot))$ then $f\in \L(X,(\clh_\centerdot))$.
\end{remark} 
\section{Improper Stochastic Integral}
Suppose $f$ is a $(\clg_\centerdot)$- predictable process such that 
\begin{equation}\label{d3}
(\clg_\centerdot)\text{-}\!\!\int f\Ind_{\{|f|\le a_n\}}\,dX\text{ is Cauchy in {\em ucp} topology whenever }a_n\uparrow\infty.\end{equation} This of course is true if $f\in \L(X, (\clg_\centerdot))$. However, it can be seen that in the example given above in Remark \ref{rm2},  while $f\not\in \L(X, (\clg_\centerdot))$, \eqref{d3} holds.

The interlacing argument yields that if for an $f$, \eqref{d3} holds, the limit in {\em ucp} of 
\[(\clg_\centerdot)\text{-}\!\!\int f\Ind_{\{|f|\le a_n\}}\,dX\]
does not depend upon the sequence $\{a_n\uparrow\infty\}$. So let $\widetilde{\L}(X, (\clg_\centerdot))$ denote the class of $(\clg_\centerdot)$- predictable processes satisfying \eqref{d3}.

For $f\in\widetilde{\L}(X, (\clg_\centerdot))$, we define the  improper integral of $f$ w.r.t.\! $X$, denoted by $(\clg_\centerdot)\text{-}\!\!\widetilde{\int} f\,dX$ as \vspace{-8mm}

\[(\clg_\centerdot)\text{-}\!\!\widetilde{\int} f\,dX=\lim_{n\rightarrow\infty}(\clg_\centerdot)\text{-}\!\!\int f\Ind_{\{|f|\le n\}}\,dX.\]

With this we have:
\vskip 2mm

\begin{proposition}\label{pr22}
Let $\clf_t\subseteq \clg_t$ for all $t$. 
 \begin{enumerate}[(i)]
\item If $f\in\L(X, (\clf_\centerdot))$ then $f\in\widetilde{\L}(X, (\clg_\centerdot))$ and \begin{equation}\label{qq30}
(\clg_\centerdot)\text{-}\widetilde{\int} f\,dX=(\clf_\centerdot)\text{-}\!\!\int f\,dX.\end{equation}
\item If $f\in\widetilde{\L}(X, (\clf_\centerdot))\cup\widetilde{\L}(X, (\clg_\centerdot))$, then  $f\in\widetilde{\L}(X, (\clf_\centerdot))\cap\widetilde{\L}(X, (\clg_\centerdot))$ and then, 
\begin{equation}\label{qq31}
(\clg_\centerdot)\text{-}\widetilde{\int} f\,dX=(\clf_\centerdot)\text{-}\widetilde{\int} f\,dX.\end{equation}
\end{enumerate}
\end{proposition}
\begin{proof}Both the parts are consequences of the observation that when the filtrations are nested, {\em i.e.} when  \eqref{z1x} holds, then the stochastic integrals w.r.t.\! the two filtrations agree for bounded $f\in\W$ and then approximating a general $f\in\W$ by $f^n=f\Ind_{\{|f|\le n\}}$ and using the definition of the improper integral, it follows that \eqref{qq30} in case $(i)$ and \eqref{qq31} in case $(ii)$ holds.
\end{proof}

With this, we have shown that when the filtrations are nested, the stochastic integrals agree if we include improper stochastic integrals.

\section{General case}
We now consider the general case of two non-comparable filtrations, {\em i.e.} when  \eqref{z1x} may not hold. 
In the rest of the section, we fix an r.c.l.l.\! process $X$ that is a $(\clf_\centerdot)$- semimartingale as well as a $(\clg_\centerdot)$- semimartingale and a process  $f$ that is bounded $(\clf_\centerdot)$- predictable as well as $(\clg_\centerdot)$- predictable. Here too, we can focus on  $f\in\W_b(\clf_\centerdot)\cap\W_b(\clf_\centerdot) $ since we have the following observation, whose proof follows  that of Proposition \ref{pr22}.

\vskip 3mm 

\begin{proposition}\label{pr3}
Let $(\clf_\centerdot)$, $(\clg_\centerdot)$, and $X$ be such that for all $f\in\W_b(\clf_\centerdot)\cap\W_b(\clg_\centerdot)$
 \begin{equation}\label{qq5}
(\clf_\centerdot)\text{-}\!\!\int f\,dX=(\clg_\centerdot)\text{-}\!\!\int f\,dX.
\end{equation} 
 \begin{enumerate}[(i)]
 \item If $f\in\L(X, (\clf_\centerdot))\cap\W(\clg_\centerdot)$, then $f\in\widetilde{\L}(X, (\clg_\centerdot))$ and
 \begin{equation}\label{qq40}
(\clg_\centerdot)\text{-}\widetilde{\int} f\,dX=(\clf_\centerdot)\text{-}\!\!\int f\,dX.\end{equation}
\item Let $f\in \W(\clf_\centerdot)\cap\W(\clg_\centerdot)$. If $f\in\widetilde{\L}(X, (\clf_\centerdot))\cup\widetilde{\L}(X, (\clg_\centerdot))$  then  $f\in\widetilde{\L}(X, (\clf_\centerdot))\cap\widetilde{\L}(X, (\clg_\centerdot))$ and  
\begin{equation}\label{qq41}
(\clg_\centerdot)\text{-}\widetilde{\int} f\,dX=(\clf_\centerdot)\text{-}\widetilde{\int} f\,dX.\end{equation}
\end{enumerate}
\end{proposition}

Let $\clh_t=\clf_t\cap\clg_t$. It follows that $X$ is $(\clh_\centerdot)$- adapted and as a consequence, $X$ is an $(\clh_\centerdot)$- semimartingale.  This follows from  Stricker's Theorem (see \cite{rlkbv2}  Theorem 4.13).  \vskip 2mm

\begin{theorem}\label{}
Suppose that 
\begin{equation}\label{d0}
f \text{ is a  bounded }(\clh_\centerdot)\text{- predictable process}\end{equation}
then 
\begin{equation}\label{d1}
(\clf_\centerdot)\text{-}\!\!\int f\,dX=(\clg_\centerdot)\text{-}\!\!\int f\,dX.\end{equation}
\end{theorem}
\begin{proof}
The proof follows from Proposition \ref{pr22} since each side in \eqref{d1} equals $(\clh_\centerdot)\text{-}\!\int f\,dX$.
\end{proof}

As a consequence, if $f$ is a left continuous process with $f\in\W_b(\clf_\centerdot)\cap \W_b(\clg_\centerdot)$,   then \eqref{d1} holds. 

However, $f$ being $(\clf_\centerdot)$- predictable and $(\clg_\centerdot)$- predictable may not imply that $f$ is $(\clh_\centerdot)$- predictable. While we do not have a counterexample to show that such an $f$ may not be $(\clh_\centerdot)$- predictable, we have not been able to prove  that $f$ is $(\clh_\centerdot)$- predictable.

Let $\clk_t=\sigma(\clf_t\cup\clg_t)$ for $t\ge 0$. If $X$ is also a $(\clk_\centerdot)\text{-}$ semimartingale, then it would follow that for all $f\in\W_b(\clf_\centerdot)\cap \W_b(\clg_\centerdot)$
\[(\clf_\centerdot)\text{-}\!\!\int f\,dX=(\clg_\centerdot)\text{-}\!\!\int f\,dX\]
since both would equal $(\clk_\centerdot)\text{-}\!\!\int f\,dX$ in view of Proposition \ref{pr22}.
\vskip 3mm

\begin{remark}
A natural question to ask is : Does $X$ being a $(\clf_\centerdot)\text{-}$ semimartingale as well as $(\clg_\centerdot)\text{-}$ semimartingale imply that $X$ is $(\clk_\centerdot)\text{-}$ semimartingale. The answer is negative as the following example shows.

Let $\{U_k:k\ge 1\}$  be a  sequence of $\{-1,1\}$ valued independent random variables with $\Pb(U_k=1)=0.5$, $\Pb(U_k=-1)=0.5$ for all $k$. 

Let $a_n=1-\frac{1}{n}$ and let $X_t$ for $t\in[0,\infty)$ be defined by
\begin{equation}\label{tr1}
X_t=\sum_{n\,:\,a_n\le t}U_{2n}U_{2n+1}\frac{1}{n}
\end{equation}
For $t<1$, $X_t$ is a sum of finitely many random variables and for $t\ge 1$, $X_t=\sum_nU_{2n}U_{2n+1}\frac{1}{n}$, a series that converges almost surely (say by Kolmogorov's 3-series theorem). Thus it follows that $X_t$ has r.c.l.l.\! paths. 
Let $\clf_t=\sigma(\{U_{2n},n\ge 1\}\cup\{X_s\,:\;s\le t\})$ and $\clg_t=\sigma(\{U_{2n+1},n\ge 1\}\cup\{X_s\,:\;s\le t\})$. It is  easy to see see that $X$ is a $(\clf_\centerdot)\text{-}$ martingale as well as a $(\clg_\centerdot)\text{-}$ martingale. In this case, $
\clk_t=\sigma(\clf_t\cup\clg_t)=\sigma(U_k\,:\;k\ge 1)=\clk_0$ for all $t\ge 0$. Thus,  all $(\clk_\centerdot)\text{-}$ martingales are constants and thus any $(\clk_\centerdot)\text{-}$ semimartingale is a process with finite variation paths. But total variation of $X_t$ over the interval $[0,1]$ equals $\sum_{n=1}^\infty \frac{1}{n}=\infty$. Thus,  $X$ is not a   $(\clk_\centerdot)\text{-}$ semimartingale.

\end{remark}

Let us return to the case of two, non comparable filtrations.
Slud \cite{Slud} had showed  that if $f\in\W_b(\clf_\centerdot)\cap\W_b(\clg_\centerdot)$ is a locally bounded process, the process $D$ defined by
\begin{equation}\label{slu}
Z_t=\Bigl((\clf_\centerdot)\text{-}\!\!\int f\,dX\Bigr)-\Bigl((\clg_\centerdot)\text{-}\!\!\int f\,dX\Bigr)
\end{equation}
is a continuous process such that $[Z,Z]=0$.

Zheng \cite{Zheng} had earlier showed that if $X$ is a semimartingale for both the filtrations and if 
\begin{equation}\label{j1}
\sum_{0\le s\le t}|(\Delta X)_s|<\infty,
\end{equation}
 then for  $f\in\W_b(\clf_\centerdot)\cap\W_b(\clg_\centerdot)$, $Z$ is a continuous process with finite variation paths. In this case the problem essentially reduced 
to that for a continuous semimartingale.

Here, we will show that under some additional conditions on  $X$, \eqref{qq0} holds.
We will now assume that \eqref{j1} holds. 
Then
\begin{equation}\label{j2}
V_t=\sum_{0\le s\le t}(\Delta X)_s<\infty
\end{equation}
is a process whose paths have finite variation and is $(\clf_\centerdot)$- adapted as well as $(\clg_\centerdot)$- adapted. Recall that $C(X,(\clf_\centerdot))$ denotes the continuous $(\clf_\centerdot)$- local martingale part of $X$.
\vskip 3mm

\begin{lemma}\label{l1}
Suppose $X$ is a $(\clf_\centerdot)$- semimartingale as well as a $(\clg_\centerdot)$- semimartingale satisfying \eqref{j1}. Let
\begin{equation}\label{ja}
D=C(X,(\clf_\centerdot))-C(X,(\clg_\centerdot)).
\end{equation}
 Then $D$ is a continuous process whose paths have finite variation. 
\end{lemma}
\begin{proof}
Let  $Y$ be defined by
\begin{equation}\label{j3}
Y_t=X_t-V_t.
\end{equation}
It follows easily that $Y$ is a continuous process. Since $V$ is a process with finite variation paths and is adapted for both the filtrations, it follows the process $Y$ is a $(\clf_\centerdot)$- semimartingale as well as a $(\clg_\centerdot)$- semimartingale.  As a consequence, $[Y,Y]$ is adapted to $(\clf_\centerdot)$ as well as $(\clg_\centerdot)$. 
Also, $[X,X]=[Y,Y]+\sum_{0<s\le t}(\Delta V)_s^2$ and hence $[Y,Y]=[X,X]^c$.

Let 
\begin{equation}\label{c2}
M=C(X,(\clf_\centerdot)),\;\;A=Y-M
\end{equation}
\begin{equation}\label{c3}
N=C(X,(\clg_\centerdot)),\;\;B=Y-N.
\end{equation}
Since $C(X,(\clf_\centerdot))=C(Y,(\clf_\centerdot))$ (as $X-Y$ is a process with finite variation paths), it follows that $Y=M+A$ is the canonical decomposition of the continuous $(\clf_\centerdot)$- semimartingale $Y$, with $M$ being a continuous $(\clf_\centerdot)$- local martingale and $A$ being $(\clf_\centerdot)$- adapted process with finite variation paths. Likewise, 
 $Y=N+B$ is the canonical decomposition of the continuous $(\clg_\centerdot)$- semimartingale $Y$, with $N$ being a continuous $(\clg_\centerdot)$- local martingale and $B$ being $(\clg_\centerdot)$- adapted process with finite variation paths. Since 
 \[D_t=M_t-N_t=B_t-A_t\]
 and $B, A$ are continuous processes with finite variation paths, it follows that $D$ also is a  continuous processes with finite variation paths.
 \end{proof}
 
 {Let us consider continuous real valued functions $F,\,G$ on $[0,\infty)$ such that $F$ is increasing, $F(0)=0$, and $G$ is a function with finite variation on $[0,T]$ for all $T<\infty$ with $G(0)=0$. We say that $G$ is absolutely continuous w.r.t. $F$ if there exists a Borel measurable function $g$ on $[0,\infty)$ such that
 \[G(t)=\int_0^tg(s)\,dF(s).\]  
}
\begin{theorem}\label{T1}
Let $X$ be a $(\clf_\centerdot)$- semimartingale as well as a $(\clg_\centerdot)$- semimartingale satisfying \eqref{j1}. Let $D=C(X,(\clf_\centerdot))-C(X,(\clg_\centerdot))$. Suppose that 
\begin{equation}\label{j4}
\Pb(\omega\,:\,s\mapsto D_s(\omega)\text{ is absolutely continous w.r.t. }s\mapsto [X,X]^c_s(\omega))=1.\end{equation}
Then for $f\in \L(X,(\clf_\centerdot))\cap \L(X,(\clg_\centerdot))$, we have
\begin{equation}\label{qq9}
(\clf_\centerdot)\text{-}\!\!\int f\,dX=(\clg_\centerdot)\text{-}\!\!\int f\,dX.
\end{equation} 
\end{theorem}
\begin{proof}In view of Proposition \ref{pr3}, 
suffices to prove that \eqref{qq9} is true for $f\in\W_b(\clf_\centerdot)\cap\W_b(\clg_\centerdot)$. So let us fix $f\in\W_b(\clf_\centerdot)\cap\W_b(\clg_\centerdot)$ and let $K$ be a bound for $f$. 

Let $V,Y$ be defined by \eqref{j2}, \eqref{j3}. Since $V$ has finite variation paths, in order to prove \eqref{qq9}, suffices to prove 
\begin{equation}\label{qz1}
(\clf_\centerdot)\text{-}\!\!\int f\,dY=(\clg_\centerdot)\text{-}\!\!\int f\,dY.
\end{equation} 
Let $M,A,N,B$ be defined by \eqref{c2}, \eqref{c3}. Since $A,B,D$ are processes with finite variation paths and the integral w.r.t.\! these processes does not depend upon the underlying filtration. Thus \eqref{qz1} is same as
\begin{equation}\label{qz2}
(\clf_\centerdot)\text{-}\!\!\int f\,dM + \int f\, dA=(\clg_\centerdot)\text{-}\!\!\int f\,dN+ \int f\, dB.
\end{equation} 
We will first prove that there exist continuous processes $\{f^n:\,n\ge 1\}$, uniformly bounded by a constant $K$ such that $f^n$ is $(\clf_\centerdot)$- adapted as well as $(\clg_\centerdot)$- adapted for all $n$ and 
\begin{equation}\label{c9}
\lim_{n\rightarrow \infty}\int_0^T|f^n_t-f_t|^2\,d[Y,Y]_t=0\;\; a.s. ,\;\;\forall T<\infty.
\end{equation} 
First, note that (see \eqref{c1})  
\begin{equation}\label{c5}
[X,X]^c=[Y,Y]=[M,M]=[N,N]
\end{equation} 
We will use the technique of {\em random change of time}, see Karandikar \cite{rlkthesis, rlk1, rlk5} and Karandikar-Rao \cite{rlkbv2}. Let 
\begin{equation}\label{c99}
\tca_s=\inf\{t\ge 0\,:\; (t+[Y,Y]_t)\ge s \}.\end{equation} 
Since $[Y,Y]$ is $ (\clf_\centerdot)$- as well as $ (\clg_\centerdot)$- adapted, it follows that  for all $s$, $\tca_s$ is a $ (\clf_\centerdot)$- as well as $ (\clg_\centerdot)$- stopping time and for all $\omega$, $s\mapsto \tca_s(\omega)$ is a strictly increasing function from $[0,\infty)$ onto $[0,\infty)$. Thus, $\tca=(\tca_s)$ is a $ (\clf_\centerdot)$- random time change as well as $ (\clg_\centerdot)$- random time change. See Karandikar-Rao \cite{rlkbv2}, Chapter 7 for all results on random time change used in this article. Let
\[\clh_s=\clf_{\tca_s},\;\;\clk_s=\clg_{\tca_s},\;\;0\le s<\infty .\]
Then $ (\clh_\centerdot)$ and $ (\clk_\centerdot)$ are filtrations and $h$ defined by
\[h_s=f_{\tca_s}, \;\;C_s=[Y,Y]_{\tca_s}\]
are $ (\clh_\centerdot)$- as well as $ (\clk_\centerdot)$- predictable, with $C$ being in addition, a continuous increasing process. Further, 
\begin{equation}\label{z1}
0\le C_t(\omega)-C_s(\omega)\le (t-s),\;\;\;\forall s<t\;\;\forall \omega .
\end{equation}
Let $\rho(t)=c\exp\{-\frac{1}{(1-t^2)}\}$ for $|t|\le 1$ and zero for $|t|>1$, where $c$ is chosen so that $\int_{-1}^1\rho(s)\,ds=1$.
For $n\ge 1$, let processes $h^n$ be defined by
\begin{equation}\label{z1a}
h^n_t(\omega)=n\int_{(t-\frac{1}{n})\vee 0}^t\rho(n(t-s-{\textstyle\frac{1}{n}}))h_s(\omega)ds.
\end{equation}
Then it follows that for each $n$, $h^n$ is a continuous process and is $ (\clh_\centerdot)$ as well as $ (\clk_\centerdot)$- predictable. Further, standard arguments involving convolution yield that
\begin{equation}\label{z2}
\int_0^T|h^n_t-h_t|\,dt \rightarrow 0 \;\;a.s.\;\;\text{ as }n\rightarrow \infty.
\end{equation}
See Friedman \cite{friedman} p56. In view of \eqref{z1}, we also conclude from \eqref{z2} that
\begin{equation}\label{z2a}
\int_0^T|h^n_t-h_t|\,dC_t \rightarrow 0 \;\;a.s.\;\;\text{ as }n\rightarrow \infty.
\end{equation}
Let 
\[\tcb_t(\omega)=\inf\{s\ge 0\,:\;\tca_s(\omega)\ge t \}\]
and
\[f^n_t(\omega)=h^n_{\tcb_t(\omega)}(\omega).\]
It follows that for each $n$, $f^n$ is a continuous process and is $ (\clh_\centerdot)$ as well as $ (\clk_\centerdot)$- predictable.  Using change of variable, the fact that $f^n,f$ are
bouned by $K$ and
using \eqref{z2a} we conclude
\[\begin{split}
\lim_{n\rightarrow \infty}\int_0^T|f^n_t-f_t|^2\,d[Y,Y]_t&\le 2K \lim_{n\rightarrow \infty}\int_0^T|f^n_t-f_t|\,d[Y,Y]_t\\
&=2K\lim_{n\rightarrow \infty}\int_0^{\tca_T}|h^n_s-h_s|\,dC_s\\
&=0\;\;\;a.s.
\end{split}\]

Having proven \eqref{c9}, we now observe that continuity of $f^n$ yields
\[
(\clf_\centerdot)\text{-}\!\!\int f^n\,dY=(\clg_\centerdot)\text{-}\!\!\int f^n\,dY.
\]
and since $Y=M+A=N+B$, we have $D=B-A=M-N$ and so we can conclude
\begin{equation}\label{z8}
(\clf_\centerdot)\text{-}\!\!\int f^n\,dM-(\clg_\centerdot)\text{-}\!\!\int f^n\,dN=\int f^n\,dD
\end{equation}
Now $[Y,Y]=[M,M]=[N,N]$ along with \eqref{c9} implies
\begin{equation}\label{z9}
(\clf_\centerdot)\text{-}\!\!\int f^n\,dM\rightarrow (\clf_\centerdot)\text{-}\!\!\int f\,dM
\end{equation}
\begin{equation}\label{z10}
(\clg_\centerdot)\text{-}\!\!\int f^n\,dN\rightarrow (\clg_\centerdot)\text{-}\!\!\int f\,dN
\end{equation}
and the assumption \eqref{j4} that $D$ is absolutely continuous w.r.t.\! $[Y,Y]$ and the fact that $f^n$, $f$ are bounded by $K$ along with \eqref{c9} implies that 
\[\int_0^T |f^n-f|\,dD\rightarrow 0\;\;a.s.,\;\forall T<\infty\]
and as a consequence 
\begin{equation}\label{z11}
\int f^n\,dD\rightarrow \int f\,dD.
\end{equation}
Now, \eqref{z8}-\eqref{z11} yield
\begin{equation}\label{z12}
(\clf_\centerdot)\text{-}\!\!\int f\,dM-(\clg_\centerdot)\text{-}\!\!\int f\,dN=\int f\,dD
\end{equation}
from which we conclude that  \eqref{qz2} holds. This completes the proof as noted earlier.
\end{proof}

\begin{remark}
For $t\ge 0$, let $\clh_t=\clf_t\cap\clg_t$ and let $A$ be a continuous $(\clh_\centerdot)$ adpated strictly increasing process with $A_0(\omega)=0$ for all $\omega$. Let $ \widetilde{\Omega}=[0,\infty)\times\Omega$ and $\widetilde{\clf}=\clb_{[0,\infty)}\otimes\clf$, where $\clb_{[0,\infty)}$ is the Borel $\sigma$-field on $[0,\infty)$. 
For $C\in \widetilde{\clf}$, let
\[\mu(C)=\int_0^\infty\int\Ind_C(t,\omega)\,dA_t(\omega)\,d\Pb(\omega).\]
Then
 \[\clp(\clf_\centerdot)^\mu\cap\clp_(\clg_\centerdot)^\mu=\clp_(\clh_\centerdot)^\mu.\]
 where $\clp(\clf_\centerdot)^\mu$, $\clp_(\clg_\centerdot)^\mu$ and $\clp_(\clh_\centerdot)^\mu$ denote the $\mu$ completions of the respective $\sigma$-fields.
\end{remark}

\begin{remark}
It is easy to see from the proof that instead of \eqref{j4}, suffices to assume that there exists a continuous increasing process $V$ such that $V$ is $(\clf_\centerdot)$- adapted as well as $(\clg_\centerdot)$- adapted  and  
\begin{equation}\label{j4x}
D=C(X,(\clf_\centerdot))-C(X,(\clg_\centerdot)) \text{ is absolutely continuous w.r.t. }V.
\end{equation}
Just use $V$ instead of $[X,X]^c$ in the definition of the time change in \eqref{c99}.
\end{remark}

\section{Lebesgue decomposition of Increasing processes}
In this section we will deduce analogue of the  Lebesgue decomposition theorem for increasing processes and for processes with finite variation paths. This will be useful in the subsequent section. The main conclusion is that the Radon Nikodym derivate can be chosen to be predictable.

First we consider increasing processes. 
\begin{lemma}\label{ld1}
Let $A,R$ be continuous $(\clf_\centerdot)$ adapted increasing processes such that $A(0)=0$, $R(0)=0$. Then there exist $(\clf_\centerdot)$- predictable processes $\phi$ and $\Gamma\in\clp(\clf_\centerdot)$ such that for all $t\in[0,\infty)$  we have
\begin{equation}\label{eq101}
A_t(\omega)=\int_0^t\phi_s(\omega)\,dR_s(\omega)+\int_0^t\Ind_{\Gamma}(s,\omega)\,dA_s(\omega)
\end{equation}
and
\begin{equation}\label{eq1222}
\int_0^t\Ind_\Gamma(s,\omega)\,dR_s(\omega)=0
\end{equation}
\end{lemma}
\begin{proof}
Recall that $\clp(\clf_\centerdot)$ is the predictable $\sigma$-field for the filtration $(\clf_\centerdot)$ on $\widetilde{\Omega}=[0,\infty)\times\Omega$. 
Let us define possibly $\sigma$-finite measures $\mu$ and $\lambda$  on $(\widetilde{\Omega},\clp(\clf_\centerdot))$ by 
\begin{equation}\label{eq112}
\mu(E)=\int_\Omega\int_0^\infty\Ind_E(t,\omega)\,dA_t(\omega)\,d\Pb(\omega)\end{equation}
\begin{equation}\label{eq122}
\lambda(E)=\int_\Omega\int_0^\infty\Ind_E(t,\omega)\,dR_t(\omega)\,d\Pb(\omega) \end{equation}
for $E\in\clp(\clf_\centerdot)$.
The Lebesgue decomposition theorem applied to $\mu$, $\lambda$ yields a predictable process $\phi$ and a predictable set $\Gamma$ such that for any $E\in\clp(\clf_\centerdot)$
\[\mu(E)=\int_E \phi \,d\lambda +\mu(E\cap \Gamma)\]
and $\lambda(\Gamma)=0$. In turn, $\lambda(\Gamma)=0$ along with \eqref{eq122} (the definition of $\lambda$) implies that \eqref{eq1222} is true for all $t$. For the other part, 
let 
\begin{equation}\label{eq1023}
B_t(\omega)=\int_0^t\phi_s(\omega)\,dR_s(\omega)+\int_0^t\Ind_{\Gamma}(s,\omega)\,dA_s(\omega).\end{equation}
By definition, $B$ is a continuous increasing process with $B_0=0$. 
To complete the proof, we will show that $B_t=A_t$ $a.e.$ for all $t$ proving \eqref{eq101}. Note that by definition, for any $E\in\clp(\clf_\centerdot)$,
\[\begin{split}
\int_\Omega\int_0^\infty&\Ind_E(t,\omega)\,dB_t(\omega)\,d\Pb(\omega)\\
&=
\int_\Omega\int_0^\infty\Ind_E(t,\omega)\phi_t(\omega)\,dR_t(\omega)\,d\Pb(\omega)+
\int_\Omega\int_0^\infty\Ind_E(t,\omega)\Ind_{\Gamma}(t,\omega)\,dA_t(\omega)
\\
&=\int_E \phi \,d\lambda +\mu(E\cap \Gamma)\\
&=\mu(E).\end{split}
\]
Taking $E=(s,t]\times C$ for $C\in\clf_s$ with $s<t$, and using the definition \eqref{eq112} of $\mu$ it follows that
\[\Eb[\Ind_C(B_t-B_s)]=\Eb[\Ind_C(A_t-A_s)]\]
And thus $A-B$ is a martingale. Since $A-B$ is a continuous process that is difference of two increasing process and is a martingale, it follows that $A_t-B_t=0$ $a.e.$ for all $t$ and thus \eqref{eq101} holds {\em a.e.} for every $t$. 
 This is a standard result in stochastic calculus, but for an elementary proof see \cite{rlk0}, \cite{rlkthesis}. Since $A,B$ are continuous process, this yields : the set
\[\Omega_0=\{\omega\,:\;A_t(\omega)=B_t(\omega)\;\;\forall t\}\]
has $\Pb(\Omega_0)=1$. Replacing $\phi$ by $\phi\Ind_{\Omega_0}$ and $\Gamma$ by $\Gamma\cup \Omega_0^c$, it follows that \eqref{eq101} holds. Since $\clf_0$ contains all $\Pb$- null sets by assumption, it follows that $\Gamma\in\clp(\clf_\centerdot)$ and $\phi$ is $(\clf_\centerdot)$- predictable.
\end{proof}
\begin{remark}\label{rk43}
Let us note that the result given above is essentially applying the Lebesgue decomposition theorem for each $\omega$. The main assertion in Theorem \ref{ld1} is that the Radon-Nikodym derivative of the absolutely continuous part and the support of the orthogonal component can be chosen so that they are predictable. 
\end{remark}

\begin{remark}\label{rk44}
The previous result holds if instead of continuity, the increasing processes $U,R$ are assumed to be predictable and have r.c.l.l.\! paths. In the last step, we will need to use that predictable martingales with finite variation paths are constant. This is also a standard result in stochastic calculus. See \cite{kru} for an elementary proof. Also, this follows from Theorem 8.40 in \cite{rlkbv2}. 
\end{remark}
\begin{remark}\label{rk45}
It is easy to deduce from Lemma \ref{ld1} that if the increasing processes $A$ and $R$ are such that 
\[A_t(\omega)-A_s(\omega)\le R_t(\omega)-R_s(\omega)\;\;\text{ for all }\;\;0\le s\le t;\;\;\omega\in\Omega\]
then there exists a $[0,1]$-valued predictable $\phi$ such that
\[ A_t(\omega)=\int_0^t\phi_s(\omega)\,dR_s(\omega).\]
\end{remark}
We need the following elementary facts about functions that have finite variation on compact intervals. Let $H:[0,\infty)\mapsto (-\infty,\infty)$ be a function such that $H_0=0$ and 
\[|H|_t=\sup_{0=s_0<s_1<..<s_m=t}\;\;\sum_{i=0}^{m-1}|H_{s_{i+1}}-H_{s_i}|<\infty\]
where the supremum above is over all finite partitions $\{0=s_0<s_1<..<s_m=t\}$ of $[0,t]$.
Let $H^+_t=\frac{1}{2}(|H|_t+H_t)$ and $H^-_t=\frac{1}{2}(|H|_t-H_t)$. Then $H^+$ and $H^-$ are increasing processes, $H_t=H^+_t-H^-_t$ and if $G^1$ and $G^2$ are increasing processes such that $H_t=G^1_t-G^2_t$, then
\begin{equation}\label{eq1019}
H^+_t-H^+_s\le G^1_t-G^1_s,\;\;\;\; H^-_t-H^-_s\le G^2_t-G^2_s,\;\;\;\;\forall\, 0\le s\le t<\infty.\end{equation}

\begin{theorem}\label{pac}
Let $U,R$ be continuous $(\clf_\centerdot)$ adapted processes such that $U(0)=0$, $R(0)=0$, $R$ is an increasing process and the variation $|U|_t(\omega)$ of the map $s\in[0,t]\mapsto U_s(\omega)$ is finite for all $t<\infty$. Then there exist $(\clf_\centerdot)$- predictable processes $\rho$, $\xi$ and a process $V$ with finite variation paths such that  for all $t\in[0,\infty)$ and for all $\omega\in\Omega$ we have
\begin{align}\label{eq1011}
\int_0^t|\rho_s(\omega)|\,dR_s(\omega)&<\infty \\
\label{eq1022}
\int_0^t|\xi_s(\omega)|\,dR_s(\omega)&=0\end{align}
\begin{align}\label{eq1012}
U_t(\omega)&=\int_0^t\rho_s(\omega)\,dR_s(\omega)+V_t(\omega)\\
\label{eq1032}
|V|_t(\omega)&=\int_0^t\xi_s(\omega)\,dU_s(\omega).
\end{align}
Further, $\xi$ takes values in the set $\{0,1,-1\}$.
\end{theorem}
\begin{proof}
Let $U^+_t=\frac{1}{2}(|U|_t+U_t)$ and $U^-_t=\frac{1}{2}(|U|_t-U_t)$. Then $U^+$ and $U^-$ are increasing (adapted continuous) processes. As seen in Remark \ref{rk45}, we can get $[0,1]$- valued $(\clf_\centerdot)$-predictable processes $\psi$ such that (writing $S_t=|U|_t$ for notational convenience) 
\begin{equation}\label{eq143}
U^+_t(\omega)=\int_0^t\psi_s(\omega)\,dS_s(\omega)\;\;\forall t
\end{equation}
Since $U^++U^-=S$, it follows that 
\begin{equation}\label{eq144}
U^-_t(\omega)=\int_0^t(1-\psi_s(\omega))\,dS_s(\omega)\;\;\forall t.
\end{equation}
We claim that 
\begin{equation}\label{eq145}
\int \int_0^t\min\{\psi_s(\omega),1-\psi_s(\omega)\}\,dS_s(\omega)\,d\Pb(\omega)=0 \;\;\forall t.
\end{equation}
{  Note that the processes $G^1,G^2$ defined by
\[G^1_t=\int_0^t\bigl(\psi_s(\omega)-\min\{\psi_s(\omega),1-\psi_s(\omega)\} \bigr)\,dS_s(\omega)\]
and
\[G^2_t=\int_0^t\bigl((1-\psi_s(\omega))-\min\{\psi_s(\omega),1-\psi_s(\omega)\} \bigr)\,dS_s(\omega)\]
are increasing processes, with $G^1_t\le U^+_t$ and $G^1_t-G^2_t=U_t$ for all $t$.} In view of the observation made just before the statement of this theorem, it follows that $G^1_t\ge U^+_t$ for all $t$ and thus we get $G^1_t= U^+_t$ for all $t$. This implies \eqref{eq145}. Thus $\psi_s(\omega)$ is $\{0,1\}$ valued. Let $\Lambda=\{(s,\omega):\,\psi_s(\omega)=1\}$. Then $\psi_s(\omega)=\Ind_\Lambda(s,\omega)$. Note $\Lambda$ is predictable as $\psi$ is. So we conclude that for all t
\begin{equation}\label{eq148}
U^+_t(\omega)=\int_0^t\Ind_\Lambda(s,\omega)\,dS_s(\omega);\;\;\;\;\;\;\;U^-_t(\omega)=\int_0^t(1-\Ind_\Lambda(s,\omega))\,dS_s(\omega).
\end{equation}
Using $S_t=U^+_t+U^-_t$ for all $t$, we can conclude
\begin{equation}\label{eq149}
U^+_t=\int_0^t\Ind_\Lambda(s,\omega)\,dU^+_s(\omega);\;\;\;\;U^-_t=\int_0^t\Ind_\Lambda^c(s,\omega)\,dU^-_s(\omega).
\end{equation}
\begin{equation}\label{eq150}
\int_0^t\Ind_\Lambda^c(s,\omega)\,dU^+_s(\omega)=0;\;\;\;\;\int_0^t\Ind_\Lambda(s,\omega)\,dU^-_s(\omega)=0.
\end{equation}
Using Theorem \ref{ld1} for $U^+$ and $R$ and $U^-$ and $R$ respectively we can get predictable $\phi^+$, $\Gamma^+$ and $\phi^-$, $\Gamma^-$  such that
\[U^+_t(\omega)=\int_0^t\phi^+_s(\omega)\,dR_s(\omega)+\int_0^t\Ind_{\Gamma^+}(s,\omega)\,dU^+_s(\omega)\]
\[U^-_t(\omega)=\int_0^t\phi^-_s(\omega)\,dR_s(\omega)+\int_0^t\Ind_{\Gamma^-}(s,\omega)\,dU^-_s(\omega)\]
and
\begin{equation}\label{eq160}
\int_0^t\Ind_{\Gamma^+}(s,\omega)\,dR_s=0=\int_0^t\Ind_{\Gamma^-}(s,\omega)\,dR_s
\end{equation}
Writing $\xi^+=\Ind_{\Gamma^+}\Ind_\Lambda(s,\omega)$, $\xi^-=\Ind_{\Gamma^-}\Ind_\Lambda^c(\omega)$, using \eqref{eq149}-\eqref{eq150}, we conclude
\[U^+_t(\omega)=\int_0^t\phi^+_s(\omega)\,dR_s(\omega)+\int_0^t\xi^+_s(\omega)\,dU_s(\omega)\]
\[U^-_t(\omega)=\int_0^t\phi^-_s(\omega)\,dR_s(\omega)+\int_0^t\xi^-_s(\omega)\,dU_s(\omega)\]
Taking $\rho=\phi^+-\phi^-$ and  $\xi=\xi^+-\xi^-$ we can see that  \eqref{eq1012}
Is true. Now \eqref{eq1022} follows from \eqref{eq160} and \eqref{eq1032} follows from the observations \eqref{eq149}-\eqref{eq150}.
\end{proof}
\begin{corollary}\label{cor1}
Let $U,R$ be continuous $(\clf_\centerdot)$ adapted processes such that $U(0)=0$, $R(0)=0$, $R$ is an increasing process and the variation $|U|_t(\omega)$ of the map $s\in[0,t]\mapsto U_s(\omega)$ is finite for all $t<\infty$. Then the following conditions are equivalent:
\begin{description}
\item[(I)] $\forall\omega\in\Omega$ the mapping $s\mapsto U_s(\omega)\text{ is absolutely continuous w.r.t. }s\mapsto R_s(\omega),$
\item[(II)] $\exists$ $(\clf_\centerdot)$- predictable process $\rho$ such that $\forall (t,\omega)\in\widetilde{\Omega}$, 
\[\int_0^t|\rho_s(\omega)|\,dR_s(\omega)<\infty \text{ and } 
U_t(\omega)=\int_0^t\rho_s(\omega)\,dR_s(\omega)
\]
\item[(III)] For any bounded $(\clf_\centerdot)$- predictable $f$, for any $T<\infty$,
\[\int_0^Tf^2_s\,dR_s=0 \text{ implies } \int_0^Tf_s\,dU_s=0 \]
\end{description}
\end{corollary}

\begin{remark}\label{rk447} We have noted above in Remark \ref{rk44} that the Lemma \ref{ld1} is true if continuity of  $U,R$ is replaced by requiring that paths are r.c.l.l.\! and the processes are predictable. The same is true of Theorem \ref{pac} and Corollary \ref{cor1}: Apart from Lemma \ref{ld1} the only additional fact needed for this is: $|U|$ is predictable if $U$ is.  See Corollary 8.24 in \cite{rlkbv2}. 
\end{remark}
 \section{Continuous Semimartingales}
As we saw in the previous section, if a semimartingale $X$ has summable jumps, {\em i.e.} $X$ 
 satisfies \eqref{j1}, then the question about equality of integrals with two filtrations reduces to that for $Y$ defined by \eqref{j2}-\eqref{j3} and $Y$ is a continuous semimartingale.
 
In this section we will focus on Continuous semimartingales. The discussion in the previous section leads us to\\
\begin{definition} A continuous $(\clf_\centerdot)$- semimartingale $Y$ is said to be regular if in the decomposition
\[Y=M+A\]
where $M$ is a $(\clf_\centerdot)$- local martingale, $A$ is a $(\clf_\centerdot)$- adapted process with finite variation paths we have  
\begin{equation}\label{j7}
\Pb(\omega\,:\,s\mapsto A_s(\omega)\text{ is absolutely continuous w.r.t. }s\mapsto [M,M]_s(\omega))=1.\end{equation}
\end{definition}

With this definition, here is a direct consequence of Theorem \ref{T1}
\begin{theorem}\label{T5}
Let $X$ be a continuous process. Let $X$ be a regular $(\clf_\centerdot)$- semimartingale as well as a regular $(\clg_\centerdot)$- semimartingale. Then for $f\in \L(X,(\clf_\centerdot))\cap \L(X,(\clg_\centerdot))$, we have
\begin{equation}\label{qq9xx}
(\clf_\centerdot)\text{-}\!\!\int f\,dX=(\clg_\centerdot)\text{-}\!\!\int f\,dX.
\end{equation} 
\end{theorem}
\begin{proof} The result follows by observing that regularity of $X$ under the two filtrations in consideration implies that condition \eqref{j4} is satisfied.
\end{proof}

\begin{example}\label{ex56} The Ito's example (see in Example \ref{ex10}) shows that  
$(W,(\clg_\centerdot))$ is a regular semimartingale where $W$ is Brownian motion, $0<T<\infty$ is fixed and $\clg_t=\sigma(W_s:0\le s\le t,\;W_T)$.
 Indeed, using the decomposition in Example \ref{ex10}) we can conclude that  if $0\le t_1<t_2<\ldots t_n<\ldots$ is a sequence such that $t_n\uparrow\infty$, and
\[
\clg_t=\sigma(W_s:0\le s\le t,\;W_{t_j},0\le j<\infty)
\]
Then
\[
M_t=W_t-\int_0^t\psi_sds
\]
is a $(\clg_\centerdot)$- martingale where 
\[
\psi_s=\frac{W_{t_i}-W_s}{t_i-s},\;\;\;\;t_{i-1}< s\le t_i,\;i\ge 1
\]
And hence 
$(W,(\clg_\centerdot))$ is a regular semimartingale. 
\end{example}

This observation along with Theorem \ref{T5} yields the following:
\begin{example}\label{ex58} Let $0=t_0\le t_1<t_2<\ldots t_n\ldots$ and
$0=s_0\le s_1<s_2<\ldots s_n<\ldots$
be sequences such that $t_n\uparrow\infty$, $s_n\uparrow\infty$ and let
\[
\clf_t=\sigma(W_s:0\le s\le t,\;W_{s_j},0\le j<\infty)\]
\[\clg_t=\sigma(W_s:0\le s\le t,\;W_{t_j},0\le j<\infty)
\]
Then for $f\in \L(W,(\clf_\centerdot))\cap \L(W,(\clg_\centerdot))$, we have
\begin{equation}\label{qq9xyz}
(\clf_\centerdot)\text{-}\!\!\int f\,dW=(\clg_\centerdot)\text{-}\!\!\int f\,dW.
\end{equation} 
\end{example}

Before we proceed we note the following equivalent description of regular semimartinagles. This follows from the Corollary \ref{cor1}.

\begin{theorem}\label{T6}
Let $X$ be a continuous $(\clf_\centerdot)$- semimartingale. Then $X$ is regular if and only if for any bounded $(\clf_\centerdot)$- predictable process $f$ and $T<\infty$
\begin{equation}\label{qa1}
\int_0^Tf^2_s\,d[X,X]_s=0\;\; a.s.\text{ implies } \int_0^Tf_s\,dX_s=0\;\;a.s.
\end{equation}
\end{theorem}
\begin{proof} Let $X$ be a regular $(\clf_\centerdot)$- semimartingale.
Let $M=C(X,(\clf_\centerdot))$ and $A=X-M$. Let  $f$ be a bounded $(\clf_\centerdot)$- predictable process such that 
\begin{equation}\label{eqa1}
\int_0^tf^2_s\,d[X,X]_s=0\;\;\;a.s.
\end{equation}
$X$ being regular $(\clf_\centerdot)$- semimartingale, using Corollary \ref{cor1}, 
it follows that 
\begin{equation}\label{eqa2}
(\clf_\centerdot)\text{-}\!\!\int_0^tf_s\,dA_s=0\;\;\;a.s.
\end{equation}
Also, \eqref{eqa1} along with $[X,X]=[M,M]$ implies 
\begin{equation}\label{eqa3}
(\clf_\centerdot)\text{-}\!\!\int_0^tf_s\,dM_s=0\;\;\;a.s.
\end{equation}
Now \eqref{eqa2} and \eqref{eqa3} together imply $(\clf_\centerdot)\text{-}\!\!\int f\,dX=0$. 

For the converse part, as seen above, if \eqref{qa1} holds, then so does
\begin{equation}\label{qa11}
\int_0^Tf^2_s\,d[X,X]_s=0\;\; a.s.\text{ implies } \int_0^Tf_s\,dA_s=0\;\;a.s.
\end{equation}
Using Corollary \ref{cor1} it follows that if \eqref{qa11} holds, then $X$ is regular.
\end{proof}

\begin{theorem}\label{T8}
Let $X$ be a continuous process. Let $X$ be a regular $(\clf_\centerdot)$- semimartingale and let $(\clk_\centerdot)$ be a filtration such that $X$ is $(\clk_\centerdot)$ adapted and 
\begin{equation}\label{eq55}
\clk_t\subseteq \clf_t\;\;\;\forall t\ge 0.
\end{equation}
Then $X$ is a regular $(\clk_\centerdot)$- semimartingale.
\end{theorem}
\begin{proof}
By Stricker's Theorem (see \cite{rlkbv2}  Theorem 4.13), it follows that  $X$ is a $(\clk_\centerdot)$- semimartingale. Since every  $(\clk_\centerdot)$- predictable process is $(\clf_\centerdot)$- predictable, the result follows from Theorem \ref{T6}.
\end{proof}

\begin{example}\label{ex580}  Let $W$ be a Brownian motion and $0=t_0\le t_1<t_2<t_n\ldots t_n$ and
$0=s_0\le s_1<s_2<s_n\ldots s_n$
 be sequences such that $t_n\uparrow\infty$, $s_n\uparrow\infty$. Let $(\clh_\centerdot)$ and $(\clk_\centerdot)$ such that 
\[
\sigma(W_s:0\le s\le t)\subseteq \clh_t\subseteq\clf_t=\sigma(W_s:0\le s\le t,\;W_{s_j},0\le j<\infty)\]
\[\sigma(W_s:0\le s\le t)\subseteq \clk_t\subseteq\clg_t=\sigma(W_s:0\le s\le t,\;W_{t_j},0\le j<\infty)
\]

Then for $f\in \L(W,(\clh_\centerdot))\cap \L(W,(\clk_\centerdot))$, we have
\begin{equation}\label{qq9xy}
(\clh_\centerdot)\text{-}\!\!\int f\,dW=(\clk_\centerdot)\text{-}\!\!\int f\,dW.
\end{equation} 
\end{example}
\begin{example}\label{ex66} Let $W$ be a Brownian motion and for $0\le s\le t<\infty$, let 
\[C_{s, t}=\sup_{s\le u\le t} (W_u-W_s),\]
 and for $t_1<\infty$ fixed let
\[\cla_t=\sigma(W_s:0\le s\le t, W_{t_1},\,C_{0,t_1}).\] 
Mansuy and Yor \cite{yor} have shown that for a suitable $(\cla_\centerdot)$- predictable process $f$, 
\[M_t=W_t-\int_0^tf_sds\]
is an $(\cla_\centerdot)$- martingale. Using that Brownian motion has independent increments, one can deduce that for $0<t_1<t_2$ fixed, and 
\[\clb_t=\sigma(W_s:0\le s\le t, W_{t_1},\,C_{0,t_1},(W_{t_2}-W_{t_1}),\,C_{t_1,t_2} )\]  
we can get a suitable $(\clb_\centerdot)$- predictable process $g$, such that
\[M_t=W_t-\int_0^tg_sds\]
is an  $(\clb_\centerdot)$- martingale. 
 In the same way, given a sequence $0=t_0\le t_1<t_2<\ldots t_n$ such that $t_n\uparrow\infty$, we can show that for
\[\clc_t=\sigma(W_s:0\le s\le t, W_{t_j},\,C_{t_{j-1},t_j},\; j\ge 1 )\]  
we can get a suitable $(\clc_\centerdot)$- predictable process $h$, such that
\[M_t=W_t-\int_0^th_sds\]
is an  $(\clc_\centerdot)$- martingale. Let us write $S_t=C_{0,t}$ and let
\[\cld_t=\sigma(W_s:0\le s\le t, W_{t_j},\,S_{t_j},\; j\ge 1 ).\]
Using $\cld_t\subseteq \clc_t$ and Theorem \ref{T8}, we can conclude that $W$ is a regular $(\cld_\centerdot)$- semimartingale. From here, one can deduce that for sequences $0=a_0\le a_1<a_2<a_n\ldots a_n$ and
$0=b_0\le b_1<b_2<b_n\ldots b_n$
 such that $a_n\uparrow\infty$, $b_n\uparrow\infty$,  with
\[\cle_t=\sigma(W_s:0\le s\le t, W_{a_j},\,S_{b_j},\; j\ge 1 )\]
$W$ is a regular $(\cle_\centerdot)$- martingale.
\end{example}
The above example along with Theorem \ref{T8} lead us to the following result.
\begin{theorem}\label{T12}
Let $W$ be a Brownian motion with $W_0=0$ and for $0\le t<\infty$, let 
\[S_{t}=\sup_{0\le u\le t} (W_u).\]
Let $\{a_n:n\ge 0\}$, $\{b_n:n\ge 0\}$, $\{c_n:n\ge 0\}$, $\{d_n:n\ge 0\}$ be  increasing sequences, each one increasing to $\infty$. Let $(\clf_\centerdot)$ and $(\clg_\centerdot)$ be filtrations such that for all $t<\infty$
\[\sigma(W_s:0\le s\le t)\subseteq\clf_t\subseteq \sigma(W_s:0\le s\le t, W_{a_j},\,S_{b_j},\; j\ge 1 )\]
\[\sigma(W_s:0\le s\le t)\subseteq\clg_t\subseteq \sigma(W_s:0\le s\le t, W_{c_j},\,S_{d_j},\; j\ge 1 ).\]
Then for $f\in\W_b(\clf_\centerdot)\cap\W_b(\clg_\centerdot)$
\begin{equation}\label{qz9}
(\clf_\centerdot)\text{-}\!\!\int f\,dX=(\clg_\centerdot)\text{-}\!\!\int f\,dX.
\end{equation} 
\end{theorem}
\section{Open questions}
While thinking about the question as to when is  \eqref{qq9} true, we came up with what appear to be open questions:

\noindent{\bf Question 1}: Let $M$ be a local martingale w.r.t.\! a filtration $ (\clf_\centerdot)$. Let $M=M^c+M^d$ be a decomposition of $M$ with $M^c$ being a continuous $ (\clf_\centerdot)$- local martingale and $M^d$ being a purely discontinuous $ (\clf_\centerdot)$- local martingale. Let $ (\clg_\centerdot)$ be another filtration such that $M$ is an $ (\clg_\centerdot)$- semimartingale.  Can we conclude that $M^c$ and $M^d$ are $ (\clg_\centerdot)$- semimartingales?

\noindent{\bf Question 2}: Let $M$ be a continuous local martingale and  $ (\clg_\centerdot)$ be a filtration such that  $M$ is a $ (\clg_\centerdot)$- semimartingale. Can we assert that  $M$ is a regular $ (\clg_\centerdot)$- semimartingale? If not it would be nice to know
some conditions for this to be true.

To the best of our information, these questions are open and
we believe that the answer to both the questions is yes.


\begin{thebibliography}{99}
 \bibitem{Bichteler} Bichteler, K. {\em Stochastic integration and $L^p$-theory of semi martingales}, Annals of Probability, 9, (1981), 49-89.
 \bibitem{delmey}
C. Dellacherie and P.A. Meyer {\em Probabilities and
Potential}, (1978), North-Holland: Amsterdam. 
\bibitem{friedman} Friedman, A. {\em Stochastic Differential Equations and Applications, Vol. 1} (1975), Academic Press, New York,
\bibitem{Ito}
K. Ito, {\em Extensions of stochastic integrals},  Proc. Int. Symp. Stochastic Differential Equations, Kyoto, 1976, Kinokuniya, Tokyo, 95-109, 1978
\bibitem{J78} Jacod, J. {\em Calcul Stochastique et Problemes de Martingales}.
 Lecture Notes in Mathematics, 714, (1979), Springer-Verlag.
\bibitem{rlkthesis} Karandikar, R. L. {\em Pathwise stochastic calculus of continuous semimartingales}, Ph.D. Thesis, Indian 
Statistical Institute, 1981.
\bibitem{rlk0} Karandikar, R. L. {\em A remark on paths of continuous martingales. }Expositione
Mathematicae, 1, (1983), 67-69.
\bibitem{rlk1} Karandikar, R. L. {\em Pathwise solution of stochastic differential
equatios}. Sankhya A, 43, (1981), 121-132.
\bibitem{rlk5} Karandikar, R. L. {\em A. S. approximation results for multiplicative stochastic
integrals,} Seminaire de Probabilites XVI, Lecture notes in Mathematics, 
 920, (1982), 384-391, Springer-Verlag.
\bibitem{rlk3} Karandikar, R. L. {\em On Metivier-Pellaumail inequality, Emery toplogy 
and Pathwise formuale in Stochastic calculus.} Sankhya A, 51, (1989), 121-143.
\bibitem{rlk4} Karandikar, R. L. {\em On pathwise stochastic integration}, Stochastic processes and their applications, 57, (1995), 
11-18.
\bibitem{rlkbv} Karandikar, R. L. and Rao, B. V. {\em On Quadratic Variation of Martingales}. Proc. Indian Academy of Sciences, 
 124, (2014), 457-469.
\bibitem{rlkbv2} Karandikar, R. L. and Rao, B. V. {\em Introduction To Stochastic Calculus}. Springer Singapore, 2018.
\bibitem{kru} Kruglov, V.M. {\em On Natural and Predictable Processes.} Sankhya A 78,  (2016), 43–51. 
\bibitem{yor} Mansuy,R. and Yor, M. 
{\em Random Times and Enlargements of Filtrations in a Brownian Setting} Lecture Notes in Math. No. 1873, Springer, Berlin, 2006. 
\bibitem{Protter} Protter, P. {\em Stochastic Integration and Differential Equations}, Second edition. (2004), Springer-Verlag. 
\bibitem{Slud} Slud, E. V. {\em Stability of stochastic integrals under change of filtration}.
Stochastic Processes and their Application, 50, ( 1994), 221-233. 
\bibitem{Zheng} Zheng, W. A. {\em Une remarque sur une meme i.s. calculee dans deux filtrations}
Seminaire de Probabilites XVIII. Lecture Notes in Math. No. 1059, Springer. Berlin, 1982/83, pp. 172-178. 
\end{thebibliography}
\end{document}